\documentclass[10pt]{amsart}
\usepackage{amssymb}
\usepackage{amsbsy}
\usepackage{amscd}
\usepackage{amsmath}
\usepackage{amsthm}

\theoremstyle{plain}
 \newtheorem{thm}{Theorem}[section]
 \newtheorem{lem}[thm]{Lemma}
 \newtheorem{prop}[thm]{Proposition}
 \newtheorem{cor}[thm]{Corollary}
 
\theoremstyle{definition}
 \newtheorem{defn}{Definition}[section]
\theoremstyle{remark}
 \newtheorem{rem}{Remark}[section]
 \newtheorem{ex}{Example}[section]
 
\def\Bbb{\mathbb}

\def\cal{\mathcal}

\newcommand{\Ext}{\operatorname{Ext}}
\newcommand{\Hom}{\operatorname{Hom}}
\newcommand{\codim}{\operatorname{codim}}
\newcommand{\im}{\operatorname{im}}
\newcommand{\Aut}{\operatorname{Aut}}
\newcommand{\rk}{\operatorname{rk}}

\newcommand{\NS}{\operatorname{NS}}
\newcommand{\coker}{\operatorname{coker}}
\newcommand{\Pic}{\operatorname{Pic}}

\newcommand{\Quot}{\operatorname{Quot}}

\newcommand{\Syst}{\operatorname{Syst}}

\font\b=cmr10 scaled \magstep5
\def\bigzerou{\smash{\lower1.7ex\hbox{\b 0}}}

\numberwithin{equation}{section}
\begin{document}

\title{
A note on moduli of vector bundles on rational surfaces
}
\author{K\={o}ta Yoshioka}
 
\address{
Department of mathematics, Faculty of Science, Kobe University,
Kobe, 657, Japan}

\email{yoshioka@math.kobe-u.ac.jp}
 \maketitle

\section{Introduction}
Let $(X,H)$ be a pair of a smooth rational surface $X$ and an
ample divisor $H$ on $X$.
Assume that $(K_X,H)<0$.
Let $\overline{M}_H(r,c_1,\chi)$ be the moduli space of semi-stable
sheaves $E$ of $\rk(E)=r$, $c_1(E)=c_1$ and $\chi(E)=\chi$. 
To consider relations between 
moduli spaces of different invariants is an interesting problem.
If $(c_1,H)=0$ and $\chi \leq 0$, then 
Maruyama \cite{Ma:4}, \cite{Ma:3} studied such relations
and constructed a contraction map
$\phi:\overline{M}_H(r,c_1,\chi) \to \overline{M}_H(r-\chi,c_1,0)$.
Moreover he showed that the image is the
Uhlenbeck compactification
of the moduli space of $\mu$-stable vector bundles.
In particular, he gave an algebraic structure on 
Uhlenbeck compactification which was topologically constructed before.
After Maruyama's result, 
Li \cite{Li:1} constructed the birational contraction
for general cases, by using a canonical determinant line bundle,
and gave an algebraic structure on Uhlenbeck compactification.
Although Maruyama's method works only for special cases,
his construction is interesting of its own.
Let us briefly recall his construction.
Let $E$ be a semi-stable sheaf of $\rk(E)=r$, $c_1(E)=c_1$ and $\chi(E)=\chi$. 
Then $H^i(X,E)=0$ for $i=0,2$.
We consider a universal extension
\begin{equation}\label{eq:uni-ext}
 0 \to E \to F \to H^1(X,E) \otimes {\cal O}_X \to 0.
\end{equation}
Maruyama showed that $F$ is a semi-stable sheaf
of $\rk(F)=r-\chi$, $c_1(F)=c_1$ and
$\chi(F)=0$.
Then we have a map $\phi:\overline{M}_H(r,c_1,\chi) \to
\overline{M}_H(r-\chi,c_1,0)$.
He showed that $\phi$ is an immersion on the open subscheme
consistings of $\mu$-stable vector bundles and
the image of $\phi$ is the Uhlenbeck compactification.
For the proof, the rigidity of ${\cal O}_X$ is essential.
In this note, we replace ${\cal O}_X$ by other rigid
and stable vector bundles $E_0$ and show that similar results hold,
if $E_0$-twisted degree $\deg_{E_0}(E):=(c_1(E_0^{\vee} \otimes E),H)=0$.  
If $H$ is a general polarization, then
we also show that $\im \phi$ is normal (Theorem \ref{thm:normal}).

We are also motivated by our study of sheaves on K3 surfaces. 
For K3 and abelian surfaces, integral functor called 
Fourier-Mukai functor gives
an equivalence of derived categories of coherent sheaves,
and under suitable conditions, we get a  
birational correspondence of moduli spaces
(cf. \cite{Y:5}, \cite{Y:7}, \cite{Y:9}).
For rational surfaces, we can rarely expect such an equivalence
(cf. \cite{Br:2}).
For example, an analogue of Mukai's reflection \cite{Mu:4} 
(which is given by \eqref{eq:uni-ext})
may lose some information.
Indeed we get our contraction map
$\phi:\overline{M}_H(r,c_1,\chi) \to
\overline{M}_H(r-\chi,c_1,0)$. 

In section \ref{sect:deg=1}, 
we also consider the relation of different moduli spaces
in the case where $\deg_{E_0}E=1$.
Then we find some relations on (virtual) Hodge numbers 
(or Betti numbers) of moduli spaces.
If $X={\Bbb P}^2$, by using known results on
Hodge numbers (\cite{E-S:1}, \cite{Y:1}),
we calculate Hodge numbers of
some low dimensional moduli spaces.
We also determine the boundary of ample cones 
in some cases.

\section{Preliminaries}\label{sect:pre}

\subsection{Twisted stability}
Let $X$ be a smooth projective surface.
Let $K(X)$ be the Grothendieck group of $X$.
For $x \in K(X)$, we set 
\begin{equation}
\gamma(x):=(\rk x,c_1(x),\chi(x)) \in {\Bbb Z} \oplus \NS(X) \oplus
{\Bbb Z}.
\end{equation}
Then $\gamma:K(X) \to {\Bbb Z} \oplus \NS(X) \oplus
{\Bbb Z}$ is a surjective homomorphism and $\ker \gamma$ is generated
by ${\cal O}_X(D)-{\cal O}_X$ and
${\Bbb C}_P-{\Bbb C}_Q$, where $D \in \Pic^0(X)$ and
$P,Q \in X$.
For $\gamma=(r,c_1,\chi) \in {\Bbb Z} \oplus \NS(X) \oplus
{\Bbb Z}$, we set $\rk \gamma=r$, $c_1(\gamma)=c_1$ and
$\chi(\gamma)=\chi$.
$K(X)$ is equipped with a
bilinear form $\chi(\;\;,\;\;)$:
\begin{equation}
\begin{matrix}
K(X) \times K(X) &\to &{\Bbb Z}\\
(x,y) & \mapsto &\chi(x,y)
\end{matrix}
\end{equation}
It is easy to see that
\begin{lem}
$\chi(x,y)=\chi(y,x)+(K_X,c_1(y^{\vee} \otimes x))$,
$x,y \in K(X)$.
\end{lem}
$\chi(\;\;,\;\;)$
induces a bilinear form on ${\Bbb Z} \oplus \NS(X) \oplus
{\Bbb Z}$. We also denote it by $\chi(\;\;,\;\;)$:
$\chi(\gamma(x),\gamma(y))=\chi(x,y)$.

Let ${\cal M}_H(\gamma)^{\mu \text{-}ss}$ (resp.
${\cal M}_H(\gamma)^{\mu \text{-}s}$) be the moduli stack of
$\mu$-semi-stable sheaves
(resp. $\mu$-stable sheaves) $E$ such that
$\gamma(E)=\gamma \in {\Bbb Z} \oplus \NS(X) \oplus
{\Bbb Z}$.

For $G \in K(X) \otimes{\Bbb Q}$ of $\rk G>0$, we define $G$-twisted rank, 
degree, and Euler characteristic of $x \in K(X) \otimes {\Bbb Q}$ by
\begin{equation}
\begin{split}
\rk_{G}(x)&:=\rk(G^{\vee} \otimes x)\\
\deg_{G}(x)&:=(c_1(G^{\vee} \otimes x),H)\\
\chi_{G}(x)&:=\chi(G^{\vee} \otimes x).
\end{split}
\end{equation}
For $t \in {\Bbb Q}_{>0}$, we get
\begin{equation}
\frac{\deg_G(x)}{\rk_G(x)}=\frac{\deg_{tG}(x)}{\rk_{tG}(x)},\;
\frac{\chi_G(x)}{\rk_G(x)}=\frac{\chi_{tG}(x)}{\rk_{tG}(x)}.
\end{equation}

We shall define $G$-twisted stability.
\begin{defn}[\cite{Y:9}]
Let $E$ be a torsion free sheaf on $X$.
$E$ is $G$-twisted semi-stable (resp. stable) with respect to $H$, if
\begin{equation}
\frac{\chi_G(F(nH)}{\rk_G(F)} \leq \frac{\chi_{G}(E(nH))}{\rk_{G}(E)},
n \gg 0
\end{equation}
for $0 \subsetneq F \subsetneq E$
(resp. the inequality is strict).
\end{defn} 
It is easy to see that the following relations hold:
\begin{equation}
\text{$\mu$-stable} \Rightarrow \text{$G$-twisted stable}
\Rightarrow \text{$G$-twisted semi-stable}
\Rightarrow \text{$\mu$-semi-stable}.
\end{equation}
For a ${\Bbb Q}$-divisor $\alpha$, we define $\alpha$-twisted stability
as ${\cal O}_X(\alpha)$-twisted stability.
This is nothing but the twisted stability introduced by Matsuki and 
Wentworth \cite{M-W:1}.
It is easy to see that $G$-twisted stability is determined by
$\alpha=\det(G)/\rk G$.
Hence $G$-twisted stability is the same as the 
Matsuki-Wentworth stability.

\begin{defn}
For $ \gamma \in {\Bbb Z} \oplus \NS(X) \oplus
{\Bbb Z}$,
let ${\cal M}_{H}^{G}(\gamma)^{ss}$ be the moduli stack
of $G$-twisted semi-stable sheaves $E$ of
$\gamma(E)=\gamma$ and
${\cal M}_{H}^{G}(\gamma)^{s}$ the open substack consisting of $G$-twisted
stable sheaves.
For usual stability, i.e,
$G={\cal O}_X$,
we denote ${\cal M}_{H}^{{\cal O}_X}(\gamma)^{ss}$ by
${\cal M}_{H}(\gamma)^{ss}$. 
\end{defn} 

\begin{thm}[\cite{M-W:1}]
There is a coarse moduli scheme
$\overline{M}_H^G(\gamma)$ of $S$-equivalence classes of
$G$-twisted semi-stable sheaves $E$ of $\gamma(E)=\gamma$.
\end{thm}

%\begin{rem}\label{rem:beta}
%Let $c_1(G)/\rk G=a H+\beta$, $a \in {\Bbb Q}, \beta \in H^{\perp}$ be
%the orthogonal decomposition. Then the  
%twisted semi-stability only depends on $\beta$, i.e,
%${\cal M}_H^{G}(\gamma)^{ss}={\cal M}_H^{\beta}(\gamma)^{ss}$.
%\end{rem}

\section{Construction of contraction map}\label{sect:contraction}

From now on, we assume that
 $(X,H)$ is a pair of a rational surface $X$ and an
ample divisor $H$ on $X$.
Then $\gamma:K(X) \to {\Bbb Z} \oplus \NS(X) \oplus
{\Bbb Z}$ is an isomorphism.
Assume that $(K_X,H)<0$.
Let $E_0$ be a exceptional vector bundle which is stable
with respect $H$.
Let $e_0 \in K(X)$ be the class of $E_0$ in $K(X)$. 
We set $\gamma_0:=\gamma(E_0)$ and $\omega:=\gamma({\Bbb C}_P)$,
$P \in X$.
We define homomorphism
$L_{e_0},R_{e_0}:K(X) \to K(X)$ by
\begin{equation}
\begin{split}
L_{e_0}(x)&:=x-\chi(x,e_0)e_0, x \in K(X),\\
R_{e_0}(x)&:=x-\chi(e_0,x)e_0, x \in K(X).
\end{split}
\end{equation}
Then the following relation holds.
\begin{lem}
$\chi(x,R_{e_0}(y))=\chi(L_{e_0}(x),y)$ for $x,y \in K(X)$.
\end{lem}

\subsection{Existence of $\mu$-stable vector bundle}

In this subsection, we shall give a sufficient condition
for ${\cal M}_H(r\gamma_0-a \omega)^{\mu \text{-}s}$ to be non-empty. 

\begin{lem}
${\cal M}_H(r\gamma_0-a \omega)^{\mu \text{-}ss}$ is smooth of 
$\dim {\cal M}_H(r\gamma_0-a \omega)^{\mu \text{-}ss}=2ra \rk E_0-r^2$.
\end{lem}

\begin{proof}
For $E \in {\cal M}_H(r\gamma_0-a \omega)^{\mu \text{-}ss}$,
$\Ext^2(E,E) \cong \Hom(E,E(K_X))^{\vee}=0$.
Hence ${\cal M}_H(r\gamma_0-a \omega)^{\mu \text{-}ss}$ is smooth and 
$\dim {\cal M}_H(r\gamma_0-a \omega)^{\mu \text{-}ss}=\dim \Ext^1(E,E)-
\dim \Hom(E,E)=-\chi(E,E)=2ra \rk E_0-r^2$.
\end{proof}

\begin{lem}\label{lem:necessary}
If ${\cal M}_H^{E_0}(r\gamma_0-a \omega)^s \ne \emptyset$,
then $r=1$ and $a=0$, or
$a \rk E_0-r \geq 0$.
\end{lem}

\begin{proof}
Let $E$ be an element of ${\cal M}_H^{E_0}(r\gamma_0-a \omega)^s$.
Since $E$ is simple and $\Ext^2(E,E)=0$,
$1 \geq \chi(E,E)=r^2-2ra \rk E_0$.
Hence $a \geq \frac{1}{2 \rk E_0}(r-\frac{1}{r}) \geq 0$.
If $\chi(E_0,E)=r-a \rk E_0 > 0$, then
there is a non-zero homomorphism
$E_0 \to E$.
Then 
\begin{equation}
\frac{1}{\rk E_0}=\frac{\chi(E_0,E_0)}{\rk E_0} \leq
\frac{\chi(E_0,E)}{r \rk E_0}=\frac{r-a \rk E_0}{r \rk E_0}.
\end{equation}
Therefore $a=0$ and $r=1$.
\end{proof}

\begin{lem}\label{lem:ev}
Let $E$ be a $\mu$-semi-stable sheaf of $\deg_{E_0}(E)=0$.
Then $ev:\Hom(E_0,E) \otimes E_0 \to E$ is injective and
$\coker(ev)$ is $\mu$-semi-stable.
\end{lem}

\begin{proof}
We set $G:=\ker(ev)$.
Assume that $G \ne 0$.
Let $G_0$ be a $\mu$-stable locally free subsheaf of $G$ such that
$\deg_{E_0} G_0=0$.
Then we get a non-zero homomorphism
$\phi:G_0 \to E_0$.
Since $G_0$ is locally free, $\phi$ must be an isomorphism.
Hence $\Hom(E_0,G_0) \ne 0$.
On the other hand, $ev$ induces an isomorphism
$\Hom(E_0,\Hom(E_0,E) \otimes E_0) \to \Hom(E_0,E)$.
Hence $\Hom(E_0,G)=0$, which is a contradiction.
Therefore $G=0$.
We next show that $I:=\coker(ev)$ is $\mu$-semi-stable.
Assume that $I$ has a torsion submodule $T$.
Then $J:=\ker(E \to I/T)$ is a submodule of $E$ 
containing $\im (ev)$.
By the $\mu$-semi-stability of $E$, $0 \leq \deg_{E_0}(J)=\deg_{E_0}(T)$.
Hence $T$ is of dimension 0.
Since $\im(ev)$ is locally free, $J=\im(ev)$.
Thus $I$ is torsion free.
Then it is easy to see that $\coker(ev)$ is $\mu$-semi-stable.
\end{proof}

\begin{cor}
If ${\cal M}_H^{E_0}(r\gamma_0-a \omega)^{\mu \text{-}ss}\ne \emptyset$,
then $a \geq 0$.
\end{cor}

\begin{proof}
If $a <0$, then $\dim \Hom(E_0,E) >r$ for 
$E \in {\cal M}_H^{E_0}(r\gamma_0-a \omega)^{\mu \text{-}ss}$.
By Lemma \ref{lem:ev}, we get a contradiction.
\end{proof}

\begin{prop}\label{prop:exist}
${\cal M}_H^{E_0}(r\gamma_0-a \omega)^{\mu \text{-}s} \ne \emptyset$,
if $r-a \rk E_0 \leq 0$.
Moreover, there is a $\mu$-stable locally free sheaf $E$
of $\gamma(E)=r\gamma_0-a \omega$.
\end{prop}

\begin{proof}
Let $W$ be a closed substack of 
${\cal M}_H^{E_0}(r\gamma_0-a \omega)^{\mu \text{-}ss}$
such that $E$ belongs to $W$ if and only if
there is a quotient $E \to G$ 
such that $(c_1(G)/\rk G,H)=(c_1(\gamma_0)/\rk \gamma_0,H)$
but $c_1(G)/\rk G \ne c_1(\gamma_0)/\rk \gamma_0$.
Let $f:E_0^{\oplus r} \to \oplus_{i=1}^a{\Bbb C}_{x_i}$,
$x_i \in X$ be a surjective homomorphism.
Then $E:=\ker f$ is $\mu$-semi-stable and does not belong to
$W$.
Hence ${\cal M}_H^{E_0}(r\gamma_0-a \omega)^{\mu \text{-}ss} \setminus W$
is a non-empty open substack of 
${\cal M}_H^{E_0}(r\gamma_0-a \omega)^{\mu \text{-}ss}$.  
For pairs of integers $(r_1,a_1)$ and
$(r_2,a_2)$ such that $r_1,r_2>0$, $a_1,a_2 \geq 0$
and $(r_1+r_2,a_1+a_2)=(r,a)$, 
let $N(r_1,a_1;r_2,a_2)$ be the substack of 
${\cal M}_H^{E_0}(r\gamma_0-a \omega)$ consisting of $E$
which fits in an exact sequence:
\begin{equation}
 0 \to E_1 \to E \to E_2 \to 0
\end{equation}
where $E_1$ is a $\mu$-stable sheaf of 
$\gamma(E_1)=r_1 \gamma_0-a_1 \omega$ and
$E_2$ is a $\mu$-semi-stable sheaf of 
$\gamma(E_2)=r_2 \gamma_0-a_2 \omega$.
By \cite[sect. 1]{D-L:1} or \cite[Lem. 5.2]{Y:8}, 
\begin{equation}
\begin{split}
\codim N(r_1,a_1;r_2,a_2)& \geq -\chi(E_1,E_2)\\
&=(a_1r_2+a_2r_1)\rk E_0-r_1r_2.
\end{split}
\end{equation}
By Lemma \ref{lem:necessary},
$(a_1+a_2) \rk E_0-(r_1+r_2) \geq 0$.
Hence if $a_1=0$ or $a_2=0$, then
we get $(a_1r_2+a_2r_1)\rk E_0-r_1r_2
\geq 0$.
If $a_1,a_2>0$, then
by using Lemma \ref{lem:necessary} again,
we see that $(a_1r_2+a_2r_1)\rk E_0-r_1r_2\geq
a_2r_1 \rk E_0>0$.
Therefore $N(r_1,a_1;r_2,a_2)$ is a proper
substack of ${\cal M}_H(r\gamma_0-a \omega)^{\mu \text{-}ss} 
\setminus W$, which implies that
${\cal M}_H(r\gamma_0-a \omega)^{\mu \text{-}s} \ne \emptyset$.
By \cite[Thm. 0.4]{Y:1},
the locus of non-locally free sheaves is of codimension 
$r\rk E_0-1>0$ (use \eqref{eq:quot}).
Hence ${\cal M}_H(r\gamma_0-a \omega)^{\mu \text{-}s}$ 
contains a locally free sheaf. 
\end{proof}

\subsection{Universal extension and the contraction map}

We define a coherent sheaf ${\cal E}$ on $X \times X $ by
the following exact sequence
\begin{equation}\label{eq:family}
0 \to {\cal E} \to p_1^*(E_0^{\vee}) \otimes p_2^*(E_0) 
\overset{ev}{\to} 
{\cal O}_{\Delta} \to 0.
\end{equation}
Then ${\cal E}$ is $p_2$-flat and
${\cal E}_x:={\cal E}_{|\{x \} \times X}$ is a $E_0$-twisted stable sheaf
of $\gamma({\cal E}_x)=\rk(E_0)\gamma(E_0)-\omega$. 
In particular $\chi(E_0,{\cal E}_x)=0$.

\begin{lem}\label{lem:WIT1}
For a $\mu$-semi-stable sheaf $E$ of $\deg_{E_0}(E)=0$,
\begin{equation}
p_{2*}({\cal E} \otimes p_1^*(E_0))=
R^2p_{2*}({\cal E} \otimes p_1^*(E_0))=0.
\end{equation}
\end{lem}

\begin{proof}
For $E \in {\cal M}_H(\gamma)^{\mu \text{-}ss}$,
Lemma \ref{lem:ev} implies that
$ev:\Hom(E_0,E) \otimes E_0 \to E$ is
injective. Hence $p_{2*}({\cal E} \otimes p_1^*(E))=0$.
Since $(K_X,H)<0$,
$\deg_{E_0}(E(-K_X))>\deg_{E_0}(E)=0$.
Hence $\Ext^2(E_0,E)=\Hom(E(-K_X),E_0)^{\vee}=0$.
Then $R^2p_{2*}({\cal E} \otimes p_1^*(E)) \cong \Ext^2(E_0,E) \otimes E_0=0$.
\end{proof}
The following is our main theorem of this section.
\begin{thm}\label{thm:contract}
Let $e \in K(X)$ be a class such that $\rk e>0$ and
$\deg_{E_0}(e)=0$.
Then we have a morphism $\phi_{\gamma(e)}:\overline{M}_H(\gamma(e)) \to 
\overline{M}_H^{E_0}(\gamma(\hat{e}))$ sending
$E$ to the $S$-equivalence class of $R^1p_{2*}({\cal E} \otimes p_1^*(E))$
and the restriction of $\phi_{\gamma(e)}$ to $M_H(\gamma)^{\mu \text{-}s,loc}$
is an immersion, where $\hat{e}=R_{e_0}(e)$ and 
$M_H(\gamma(e))^{\mu \text{-}s,loc}$ is the open subscheme
consisting of $\mu$-stable vector bundles.
If $\oplus_i E_i$ is the $S$-equivalence class of $E$ with respect to
$\mu$-stability, then $\phi_{\gamma(e)}(E)$ is uniquely 
determined by $\oplus_i E_i^{\vee \vee}$ and the location of pinch points of
$\oplus_i E_i$. 
\end{thm}
In order to prove this theorem, we prepare some lemmas.

\begin{lem}
\begin{equation}
{\bf R} p_{2*}({\cal E} \otimes p_1^*(E_0)))=0.
\end{equation}
\end{lem}

\begin{proof}
By \eqref{eq:family}, we have an exact sequence
\begin{equation}
\Hom(E_0,E_0) \otimes E_0 \overset{ev}{\to} E_0 \to 
R^1 p_{2*}({\cal E} \otimes p_1^*(E_0))) \to \Ext^1(E_0,E_0) \otimes E_0.
\end{equation}
Since $ev$ is isomorphic and $\Ext^1(E_0,E_0)=0$,
we get that $R^1 p_{2*}({\cal E} \otimes p_1^*(E_0))=0$.
Therefore we get our claim.
\end{proof}

\begin{lem}\label{lem:hom=0}
For a $\mu$-semi-stable sheaf $E$ of $\deg_{E_0}(E)=0$,
\begin{equation}
\Hom(E_0,R^1 p_{2*}({\cal E} \otimes p_1^*(E)))=0.
\end{equation}
\end{lem}

\begin{proof}
By Leray spectral sequence and projection formula,
\begin{equation}
\Hom(E_0,R^1 p_{2*}({\cal E} \otimes p_1^*(E)))=
H^1(X \times X,{\cal E} \otimes p_1^*(E) \otimes p_2^*(E_0^{\vee})).
\end{equation}
Since ${\bf R} p_{1*}({\cal E} \otimes p_2^*(E_0^{\vee}))=0$,
${\bf R} p_{1*}({\cal E}\otimes p_1^*(E) \otimes p_2^*(E_0^{\vee}))
={\bf R} p_{1*}({\cal E} \otimes p_2^*(E_0^{\vee})) 
\overset{\bf L}{\otimes} E=0$.
\end{proof}
For simplicity, we set $\widehat{E}:=R^1 p_{2*}({\cal E} \otimes p_1^*(E))$.

\begin{prop}\label{prop:contract}
For a $\mu$-semi-stable sheaf $E$ of $\deg_{E_0}(E)=0$,
$\widehat{E}$ is a $E_0$-twisted semi-stable sheaf of
$\chi(E_0,\widehat{E})=0$.
\end{prop}

\begin{proof}
By \eqref{eq:family}, $\widehat{E}$ fits in an exact sequence
\begin{equation}
0 \to
\Hom(E_0,E) \otimes E_0 \overset{ev}{\to}  
 E \to \widehat{E} \to \Ext^1(E_0,E) \otimes E_0 \to 0
\end{equation}
By Lemma \ref{lem:ev}, $ \widehat{E}$ is $\mu$-semi-stable.
It is easy to see that $\chi(E_0,\widehat{E})=0$.
Assume that $\widehat{E}$ is not semi-stable and let $G$ be a destabilizing
subsheaf.
Then $\chi(E_0,G)/\rk G>0$.
By our assumption on $H$, $\Ext^2(E_0,G)=0$.
Hence $\Hom(E_0,G) \ne 0$, which contradicts to Lemma \ref{lem:hom=0}.
\end{proof}

\begin{rem}
If $E$ is $E_0$-twisted semi-stable such that
$\chi(E_0,E) \leq 0$, then
$\widehat{E}$ fits in an exact sequence
\begin{equation}\label{eq:univ.ext}
0 \to 
 E \to \widehat{E} \to \Ext^1(E_0,E) \otimes E_0 \to 0.
\end{equation}
By Lemma \ref{lem:hom=0}, \eqref{eq:univ.ext} is a universal extension.
\end{rem}

\begin{lem}
Let $E$ be a $\mu$-stable vector bundle of $\deg_{E_0}(E)=0$.
Then $\widehat{E}$ is $E_0$-twisted stable.
\end{lem}
 
\begin{proof}
We may assume that $E \ne E_0$.
Then $\widehat{E}$ fits in a universal extension
\begin{equation}
0 \to E \to \widehat{E} \to E_0^{\oplus h} \to 0
\end{equation}
where $h=\dim \Ext^1(E_0,E)$.
Assume that $\widehat{E}$ is not $E_0$-twisted stable.
Then there is a $E_0$-twisted stable subsheaf $G_1$ of $\widehat{E}$
such that $G_2:=\widehat{E}/G_1$ is $E_0$-twisted semi-stable.
If $E$ is contained in $G_1$, then 
we get a homomorphism $E_0^{\oplus h} \to G_2$.
Since $\chi(E_0, G_2)/\rk G_2=0<\chi(E_0,E_0^{\oplus h})/h \rk E_0$,
we get a contradiction.
Hence $E$ is not contained in $G_1$.
Since $E$ is $\mu$-stable, we get $E \cap G_1=0$.
Hence $G_1 \to E_0^{\oplus h}$ is injective.
%We take a non-zero homomorphism
%$G_1 \to E_0$.
Let $G'$ be a $\mu$-stable locally free subsheaf of $G_1$.
Then we see that $G' \cong E_0$, which implies that
$G_1$ is not $E_0$-twisted stable.
Therefore $\widehat{E}$ is $E_0$-twisted stable. 
\end{proof}

{\it Proof of Theorem \ref{thm:contract}:}
Let $\{{\cal F}_s \}_{s \in S}$ be a flat family of $\mu$-semi-stable
sheaves of $\deg_{E_0}({\cal F}_s)=0$.
Then Lemma \ref{lem:WIT1} and Proposition \ref{prop:contract}
imply that
$\{\widehat{{\cal F}_s} \}_{s \in S}$ is also a flat family of
$E_0$-twisted semi-stable sheaves (cf. \cite[Thm. 1.6]{Mu:5}).
Hence we get a morphism $\phi_{\gamma(e)}:
\overline{M}_H(\gamma(e)) \to \overline{M}_H(\gamma(\hat{e}))$.
Let $E$ be a $\mu$-stable vector bundle of $\deg_{E_0}(E)=0$ 
and $\varphi:E \to T$ be a 
quotient such that $T$ is of dimension 0.
Then for $F:=\ker \varphi$, we get an exact sequence
\begin{equation}
 0 \to p_{2*}({\cal E} \otimes p_1^*(T)) \to \widehat{F} \to
 \widehat{E} \to 0.
\end{equation}
Let 
$0 \subset T_1 \subset T_2 \subset \dots \subset
T_n=T$ be a filtration such that $T_i/T_{i-1} \cong {\Bbb C}_{x_i}$,
$x_i \in X$
(i.e, Jordan-H\"{o}lder filtration with respect to Simpson's stability).
Then $G:=p_{2*}({\cal E} \otimes p_1^*(T))$ has a filtration
$0 \subset G_1 \subset G_2 \subset \dots \subset
G_n=G$ such that 
$G_i/G_{i-1} \cong {\cal E}_{x_i}$.
Since $\widehat{E}$ is stable,
the $S$-equivalence class of $\widehat{F}$ is $\widehat{E} \oplus
\oplus_{i=1}^n {\cal E}_{x_i}$.  

For a $\mu$-semi-stable sheaf $E$ of $\deg_{E_0}(E)=0$,
let $\oplus_{i=1}^n E_i$ be an $S$-equivalence class of $E$
with respect to $\mu$-stability.
Let $\oplus_j {\Bbb C}_{x_{{i,j}}}$ be the $S$-equivalence class of 
$E_i^{\vee \vee}/E_i$ as a purely 0-dimensional sheaf.
Then the $S$-equivalence class of $\widehat{E}$
with respect to $E_0$-twisted stability is 
$\oplus_{i=1}^n(\widehat{E_i^{\vee \vee}} \oplus \oplus_{j}
{\cal E}_{x_{{i,j}}})$.
By Proposition \ref{prop:classification} and 
Remark \ref{rem:classification} below,
$\widehat{E_i^{\vee \vee}}$ is uniquely determined by
$E_i^{\vee \vee}$.
Hence the $S$-equivalence class of $\widehat{E}$ is 
uniquely determined by $E_i^{\vee \vee}$ and $x_{i,j}$.
\qed

\begin{prop}\label{prop:classification}
Let $F$ be an $E_0$-twisted stable sheaf such that
$\deg_{E_0}(F)=0$ and $\chi(E_0,E)=0$.
Then 
\begin{enumerate}
\item[(1)]
$F={\cal E}_x$, $x \in X$, or
\item[(2)]
$F$ fits in an exact sequence
\begin{equation}\label{eq:classification}
0 \to E \to F \to E_0^{\oplus n} \to 0,
\end{equation}
where $E$ is a $\mu$-stable locally free sheaf.
\end{enumerate}
\end{prop}

\begin{proof}
If $F$ is $\mu$-stable, then we see that
$F^{\vee \vee} \cong E_0$, and hence
$\rk E_0=1$ and $F \cong {\cal E}_x$, $x \in X$.
Assume that there is an exact sequence
\begin{equation}
0 \to G_1 \to F \to G_2 \to 0,
\end{equation}
where $G_1$ is a $\mu$-stable sheaf of $\deg_{E_0}(G_1)=0$
and $G_2$ is a $\mu$-semi-stable sheaf of $\deg_{E_0}(G_2)=0$.
Then we get an exact sequence
\begin{equation}
0 \to \widehat{G_1} \to \widehat{F} \to \widehat{G_2} \to 0.
\end{equation}
Since $F$ is $E_0$-twisted stable, $\widehat{F}=F$.
In particular $\widehat{F}$ is $E_0$-twisted stable. 
By the stability of $G_1$, 
$\chi(E_0,G_1)<0$, which implies that $\widehat{G_1} \ne 0$.
Therefore $\widehat{G_1} \cong \widehat{F}$ and $\widehat{G_2}=0$.
By using \eqref{eq:family}, we see that
$\Hom(E_0,G_2) \otimes E_0 \to G_2$ is an isomorphism.
We note that
$\widehat{G_1}$
fits in an exact sequence
\begin{equation}
0 \to p_{2*}({\cal E} \otimes p_1^*(G_1^{\vee \vee}/G_1)) \to
\widehat{G_1} \to  \widehat{G_1^{\vee \vee}} \to 0.
\end{equation}
By the stability of $F$, (i) $G_1^{\vee \vee}/G_1=0$, or
(ii) $G_1^{\vee \vee}/G_1={\Bbb C}_x$, $x \in X$ and
$\widehat{G_1^{\vee \vee}}=0$.
Therefore $G_1$ is locally free, or
$F={\cal E}_x$. 
\end{proof}

\begin{rem}\label{rem:classification}
If $F$ fits in the exact sequence
\eqref{eq:classification}, then
$E=\ker(F \to \Hom(F,E_0)^{\vee} \otimes E_0)$.
Thus $E$ is uniquely determined by $F$.
\end{rem}

\vspace{1pc}

\begin{ex}
Assume that $(X,H)=({\Bbb P}^2,{\cal O}_{{\Bbb P}^2}(1))$ and 
$E_0=\Omega_X(1)$.
Then we have a contraction
\begin{equation}
\overline{M}_H(2,-H,-n) \to \coprod_{0 \leq k \leq n}
M_H(2,-H,-k)^{\mu \text{-}s,loc} \times S^{n-k} X
\end{equation}
sending $E$ to $(E^{\vee \vee},gr(E^{\vee \vee}/E))$,
where $gr(E^{\vee \vee}/E)$ is the $S$-equivalence class of $E^{\vee \vee}/E$.
\end{ex}

\begin{rem}
For a $\mu$-semi-stable sheaf $E$ of $\deg_{E_0}(E)=0$,
${\cal H}(E):=\Ext^1_{p_1}(p_2^*(E),{\cal E})$ is a semi-stable sheaf 
such that $\deg_{E_0^{\vee}}({\cal H}(E))=0$ and
$\chi(E_0^{\vee},{\cal H}(E))=0$.
Indeed, it is easy to see that ${\cal H}(E)$ is a $\mu$-semi-stable
sheaf such that $\deg_{E_0^{\vee}} {\cal H}(E)=0$
and $\chi(E_0^{\vee},{\cal H}(E))=0$.
Since $\Hom(E_0^{\vee},{\cal H}(E))=
\Ext^1(p_2^*(E),{\cal E} \otimes p_1^*(E_0))=0$,
${\cal H}(E)$ is semi-stable. 
Hence we have a morphism
$\psi_{\gamma}:
\overline{M}_H^{E_0}(\gamma(e)) \to
\overline{M}_H^{E_0^{\vee}}(\gamma(\hat{e}^{\vee}))$.
It is easy to see that
$\psi_{\delta}$ is an isomorphism and we get a commutative
diagram.
\begin{equation}
\begin{matrix}
           && \overline{M}_H^{E_0}(\gamma(e))&&&&&
    \overline{M}_{H}^{E_0^{\vee}}(\gamma(e^{\vee}))\cr
          &\llap{$\phi_{\gamma(e)}$}
    \swarrow&&\searrow\rlap{$\psi_{\gamma(e)}$}&&  
     \llap{$\phi_{\gamma(e^{\vee})}$}\swarrow &\cr
          \overline{M}_H^{E_0}(\gamma(\hat{e}))&& 
    \overset{\psi_{\gamma(\hat{e})}}{\to}
   &&\overline{M}_H^{E_0^{\vee}}(\gamma(\hat{e}^{\vee}))&&\cr    
\end{matrix}
\end{equation}

\end{rem}

\section{The image of the contraction}\label{sect:image}

\subsection{Brill-Noether locus}

We set $\widehat{\gamma}:=m \gamma_0-c \omega$.
Assume that $H$ is general with respect to $\widehat{\gamma}$,
that is, $H$ does not lie on walls with respect to $\widehat{\gamma}$
(cf. \cite{M-W:1}, \cite{Y:2},\cite{Y:8}).
Hence ${\cal M}_H^{E_0}(\widehat{\gamma})^{ss}={\cal M}_H(\widehat{\gamma})^{ss}$.
We define Brill-Noether locus by
\begin{equation}
 {\cal M}_H(\widehat{\gamma},n):=
 \{F \in {\cal M}_H(\widehat{\gamma})^{ss}|\dim \Hom(F,E_0) \geq n \}
\end{equation}
and the open substack 
${\cal M}_H(\widehat{\gamma},n)_0={\cal M}_H(\widehat{\gamma},n) \setminus{\cal M}_H(\widehat{\gamma},n+1)$.
By using determinantal ideal,
${\cal M}_H(\widehat{\gamma},n)$ has a substack structure.
Indeed, let $Q(\widehat{\gamma})$ be a standard open covering of
${\cal M}_H(\widehat{\gamma})^{\mu \text{-}ss}$,
that is, $Q(\widehat{\gamma})$ is an open subscheme of a quot-scheme
$\Quot_{{\cal O}_X(-k)^{\oplus N}/X/{\Bbb C}}$,
$k \gg 0$, 
$N=\chi(\widehat{\gamma}(k))$
whose points consist of quotients ${\cal O}_X(-k)^{\oplus N} \to F$
such that 
\begin{enumerate}
\item
$F \in {\cal M}_H(\widehat{\gamma})^{\mu \text{-}ss}$,
\item
$H^0(X,{\cal O}_X^{\oplus N}) \to H^0(X,F(k))$ is an isomorphism
and $H^i(X,F(k))=0$ for $i>0$.
\end{enumerate}
We may assume that
\begin{equation}\label{eq:k>>0}
 H^i(X,E_0(k))=0, i>0.
\end{equation}
Let ${\cal O}_{Q(\widehat{\gamma}) \times X}(-k)^{\oplus N} \to {\cal Q}$ 
be the universal quotient and
${\cal K}$ the universal subsheaf.
We set 
\begin{equation}
\begin{split}
V:&=
\Hom_{p_{Q(\widehat{\gamma})}}({\cal O}_{Q(\widehat{\gamma}) \times X}(-k)^{\oplus N},
{\cal O}_{Q(\widehat{\gamma})} \otimes E_0),\\ 
W:&=\Hom_{p_{Q(\widehat{\gamma})}}({\cal K},{\cal O}_{Q(\widehat{\gamma})} \otimes E_0).
\end{split}
\end{equation}
Since $\Ext^2({\cal Q}_q,E_0)=0$ for all $i>0$ and $q \in Q(\widehat{\gamma})$,
\eqref{eq:k>>0} implies that $\Ext^i({\cal K}_q,E_0)=0$ for all 
$q \in Q(\widehat{\gamma})$.
Hence $V$ and $W$  
are locally free sheaves on $Q(\widehat{\gamma})$ and
we have an exact sequence
\begin{equation}
 0 \to \Hom({\cal Q}_q,E_0) \to
 V_q \to W_q \to \Ext^1({\cal Q}_q,E_0) \to 0, q \in Q(\widehat{\gamma}).
\end{equation}
Therefore we shall define the stack structure on
${\cal M}_H(\widehat{\gamma},n)$ as the zero locus of
$\wedge^r V \to \wedge^r W$.

Let ${\cal M}_H(\widehat{\gamma},n\gamma_0)$ be the moduli stack of
isomorphism classes of 
$F \to E_0^{\oplus n}$ such that
$F \in {\cal M}_H(\widehat{\gamma})^{\mu \text{-}ss}$ and
$\Hom(E_0^{\oplus n},E_0) \to \Hom(F,E_0)$ is injective.
We have a natural projection
${\cal M}_H(\widehat{\gamma},n\gamma_0) \to {\cal M}_H(\widehat{\gamma},n)$.
Let ${\cal M}_H(\widehat{\gamma},n\gamma_0)_0$ be the open substack of
${\cal M}_H(\widehat{\gamma},n\gamma_0)$ such that
$\Hom(E_0^{\oplus n},E_0) \to \Hom(F,E_0)$ is isomorphic.
By \cite[Chap.II sect. 2,3]{ACGH:1}, 
${\cal M}_H(\widehat{\gamma},n\gamma_0)_0$ is isomorphic to
${\cal M}_H(\widehat{\gamma},n)_0$.

We shall show that ${\cal M}_H(\widehat{\gamma},n)$ is Cohen-Macaulay and normal.
By \cite[Chap.II Prop.(4.1)]{ACGH:1},
if ${\cal M}_H(\widehat{\gamma},n)$ has an expected codimension, that is, 
$\codim {\cal M}_H(\widehat{\gamma},n) =n^2$, then 
${\cal M}_H(\widehat{\gamma},n)$ is Cohen-Macaulay.
We shall estimate the dimension of substack ${\cal M}_H(\widehat{\gamma};n,p,a)$ of 
${\cal M}_H(\widehat{\gamma})^{\mu \text{-}ss}$ consisting of
$F \in {\cal M}_H(\widehat{\gamma})^{\mu \text{-}ss}$ such that
$\dim F^{\vee \vee}/F=p$ and
$F^{\vee \vee}$ fits in an exact sequence
\begin{equation}
0 \to E \to F^{\vee \vee} \to G \to 0
\end{equation}
where $E$ is a $\mu$-semi-stable sheaf of $\gamma(E)=
r\gamma_0-b\omega$, $G^{\vee \vee} \cong E_0^{\oplus n}$ and 
$\gamma(G)=n\gamma_0-a\omega$.  
\begin{lem}\label{lem:estimate}
 $\codim {\cal M}_H(\widehat{\gamma};n,p,a) \geq n^2+(r \rk E_0-1)(a+p)$.
\end{lem}

\begin{proof}
For a locally free sheaf $L$,
\cite[Thm. 0.4]{Y:1} implies that
\begin{equation}\label{eq:quot}
\dim \Quot_{L/X/{\Bbb C}}^a=(\rk L+1)a.
\end{equation}
Let $N$ be the substack of 
${\cal M}_H((r+n)\gamma_0-(a+b)\omega)^{\mu \text{-}ss}$
consisting of $F$ which fits in an exact sequence
\begin{equation}
0 \to E \to L \to G \to 0
\end{equation}
where $E$ is a $\mu$-semi-stable sheaf of $\gamma(E)=
r\gamma_0-b\omega$, $G^{\vee \vee} \cong E_0^{\oplus n}$ and
$\gamma(G)=n\gamma_0-a\omega$.  
By \cite[Lem. 5.2]{Y:8}, we see that
\begin{equation}
\begin{split}
\dim N& \leq \dim {\cal M}_H(r \gamma_0-b \omega)^{\mu \text{-}ss}+
\dim ([\Quot_{E_0^{\oplus n}/X/{\Bbb C}}^a/\Aut(E_0^{\oplus n})])-
\chi(G,E)\\
&=
(2rb \rk E_0-r^2)+((n \rk E_0+1)a-n^2)+
((ra+nb)\rk E_0-rn))\\
&=(r+n)((a+b)\rk E_0-(r+n))+n(r+n)+a+br \rk E_0-n^2.
\end{split}
\end{equation}
Hence by using \eqref{eq:quot} and the assumption
$(a+b+p)\rk E_0=r+n$, we see that
\begin{equation}
\begin{split}
\dim {\cal M}_H(\widehat{\gamma};n,p,a)&=\dim N+((r+n)\rk E_0+1)p\\
& \leq n(r+n)+a+p+br \rk E_0-n^2.
\end{split}
\end{equation}
Therefore we get
\begin{equation}
\begin{split}
\codim {\cal M}_H(\widehat{\gamma};n,p,a) &\geq (r+n)(2(a+b+p) \rk E_0-(r+n))
-(n(r+n)+a+p+br \rk E_0-n^2)\\
&=(r+n)^2-n(r+n)-(a+p+br \rk E_0-n^2)\\
&=n^2+(r \rk E_0-1)(a+p).
\end{split}
\end{equation}
\end{proof}

\begin{cor}
If $r:=m-n \geq 1$, then ${\cal M}_H(\widehat{\gamma};n)$ is Cohen-Macaulay.
\end{cor}

%\begin{proof}
%If $m=n$, then
%$\codim {\cal M}_H(\widehat{\gamma};m,a) \geq m^2-a$,
%where $a \rk E_0=m$.
%Since $\codim {\cal M}_H(\widehat{\gamma};m,a)-(m-1)^2
%\geq m^2-a-(m-1)^2=2m-1-a=(m-1)+(\rk E_0-1)a \geq 1$.
%Hence $\codim {\cal M}_H(\widehat{\gamma};n) \geq n^2$ for
%$n<m$. Therefore ${\cal M}_H(\widehat{\gamma};n)$ is Cohen-Macaulay.
%\end{proof}

Assume that $r \rk E_0 \geq 2$.
Since $\codim_{{\cal M}_H(\widehat{\gamma};n)} {\cal M}_H(\widehat{\gamma};n+1)
\geq 2n+1$, we shall show that
${\cal M}_H(\widehat{\gamma};n)_0 \cong {\cal M}_H(\widehat{\gamma},n\gamma_0)_0$ 
is regular in codimension 1.
For an element $ F \to E_0^{\oplus n}$ of
${\cal M}_H(\widehat{\gamma},n\gamma_0)_0$,
the obstruction for smoothness belongs to
$\Ext^2(F,F \to E_0^{\oplus n})$.

\begin{lem}\label{lem:obstruction}
If $F \to E_0^{\oplus n}$ is surjective or 
$F$ is locally free, then
$\Ext^2(F,F \to E_0^{\oplus n})=0$.
\end{lem}

\begin{proof}
We have an exact sequence
\begin{equation}
\Ext^2(F,E) \to
\Ext^2(F,F \to E_0^{\oplus n}) \to
\Ext^2(F,G \to E_0^{\oplus n}),
\end{equation}
where $G:=\im(F \to E_0^{\oplus n})$.
Then $\Ext^2(F,E)=\Hom(E,F(K_X))^{\vee}=0$.
Since $\Ext^2(F,G \to E_0^{\oplus n})=
\Ext^1(F,E_0^{\oplus n}/G)$,
we get our claim.
\end{proof}

If $a+p \geq 2$, then 
$\codim {\cal M}_H(\widehat{\gamma};n,p,a) \geq 2$.
If $a+p \leq 1$, then 
Lemma \ref{lem:obstruction} implies that  
${\cal M}_H(\widehat{\gamma},n)$ is smooth on 
${\cal M}_H(\widehat{\gamma};n,p,a)$.
Hence ${\cal M}_H(\widehat{\gamma},n)$ is regular in codimension 1.
By Serre's criterion, ${\cal M}_H(\widehat{\gamma};n)$ is normal.

\begin{prop}\label{prop:image}
Assume that $r \rk E_0 \geq 2$. Then
${\cal M}_H(\widehat{\gamma};n)$, 
$n:=m-r$ is normal and general member $F$ fits in
an exact sequence
\begin{equation}
 0 \to E \to F \to E_0^{\oplus n} \to 0,
\end{equation}
where $E \in {\cal M}_H(\widehat{\gamma}-n\gamma_0)^{\mu \text{-}s, loc}$ and
$\Hom(E,E_0)=0$.
\end{prop}
The following is a partial answer to \cite[Question 6.5]{Ma:3}.  
\begin{thm}\label{thm:normal}
Assume that $r \rk E_0 \geq 2$.
For $n:=m-r$,
we set
\begin{equation}
\overline{M}_H(\widehat{\gamma};n):=
\{F \in \overline{M}_H(\widehat{\gamma})|\dim \Hom(F,E_0) \geq n \}.
\end{equation}
Then 
$\overline{M}_H(\widehat{\gamma};n)$ is normal, 
$\overline{M}_H(\widehat{\gamma};n)=\phi_{\gamma}(\overline{M}_H(\gamma))$ and
we have an identification
\begin{equation}
\overline{M}_H(\widehat{\gamma};n)
=\coprod_{r_i,a_i,n_i,l}
\prod_i S^{n_i} M_H(r_i \gamma_0-a_i \omega)^{\mu \text{-}s, loc} 
\times S^l X
\end{equation}
where 
$r_i,a_i,n_i,l$ satisfy that
$a_i \rk E_0 \geq r_i$,
$(r_i,a_i)  \ne (r_j,a_j)$ for $i \ne j$,
$l+\sum_i n_i a_i=a$ and
$\sum_i n_i r_i \leq r=m-n$
Therefore $\phi_{\gamma}(\overline{M}_H(\gamma))$ is normal.
\end{thm}
 
\begin{proof}
By Proposition \ref{prop:image},
$\overline{M}_H(\widehat{\gamma};n)$ is normal.
Moreover $\phi_{\gamma}(M_H(\gamma)^{\mu \text{-}s, loc})$
is a dense subset of $\overline{M}_H(\widehat{\gamma};n)$.
Hence $\overline{M}_H(\widehat{\gamma};n)=
\phi_{\gamma}(\overline{M}_H(\gamma))$.
Let $F$ be a poly-stable sheaf of $\gamma(F)=\widehat{\gamma}$,
i.e, $F$ is a direct sum of $E_0$-twisted stable sheaves.
By Proposition \ref{prop:classification},
there are $\mu$-stable locally free sheaves $E_i$, $1 \leq i \leq k$ of
$\gamma(E_i)=r_i \gamma_0-a_i \omega$ and points $x_j \in X$, $1 \leq j \leq
l$ such that  
$F=\oplus_{i=1}^k\widehat{E_i} \oplus \oplus_{j=1}^l {\cal E}_{x_j}$.
Since $\dim \Hom (\widehat{E_i},E_0)=a_i \rk E_0-r_i$ and
$\dim \Hom( {\cal E}_{x_j},E_0)=\rk E_0$, we see that
\begin{equation}
\begin{split}
\dim \Hom(F,E_0)&=\sum_i(a_i \rk E_0-r_i)+l \rk E_0\\
&=a \rk E_0-\sum_i r_i=m-\sum_i r_i.
\end{split}
\end{equation}
Hence  
$F$ belongs to $\overline{M}_H(\widehat{\gamma};n)$ if and only if
$\sum_i r_i \leq m-n=r$.
Then the last claim follows from this. 
\end{proof}

\section{The case where $\deg_{E_0}(E)=1$}\label{sect:deg=1}

\subsection{Twisted coherent systems and correspondences} 

In this section, we shall treat the case where 
the twisted degree is $1$. This case was highly motivated by
Ellingsrud and Str\o mme's paper \cite{E-S:1}. 
Assume that $\rk e_0 (-K_X,H)>1$.
Let $e$ be a class in $K(X)$ such that
$\rk e>0$ and $\deg_{e_0}(e)=1$.
We set $\gamma:=\gamma(e)$ and $\gamma_0:=\gamma(e_0)$.
For a stable sheaf $E$ of $\gamma(E)=\gamma$,
$\Hom(E,E_0)=0$.
Since $\deg_{E_0}(E(K_X))=\deg_{E_0}(E)+\rk E \rk E_0 (K_X,H)<0$,
we get $\Ext^2(E,E_0)=\Hom(E_0,E(K_X))^{\vee}=0$.
Hence $-\chi(e,e_0) \geq 0$.

\begin{prop}\label{prop:fine}
$M_H(\gamma)$ is compact and
there is a universal family on
$M_H(\gamma) \times X$.
\end{prop}

\begin{proof}
Since $\deg_{e_0}(e)=\rk e_0 (c_1(e),H)-
\rk e (c_1(e_0),H)=1$,
$\rk e$ and $(c_1(e),H)$ are relatively prime.
Hence there is a universal family.
\end{proof}
In order to construct a correspondence, we consider $E_0$-twisted
coherent systems.
Let $\Syst(E_0^{\oplus n},\gamma)$ be the moduli space of
$E_0$-twisted coherent systems:
\begin{equation}
\Syst(E_0^{\oplus n},\gamma):=\{(E,V)|E \in M_H(\gamma), V \subset
\Hom(E_0,E), \dim V=n \}.
\end{equation} 
$\Syst(E_0^{\oplus n},\gamma)$ is a projective scheme
over $M_H(\gamma)$ (cf. \cite{L:1}).

We set 
\begin{equation}
M_H(\gamma)_i:=\{E \in M_H(\gamma)|
\dim \Hom(E_0,E)=i \}.
\end{equation}
If $i \geq n$, then $\Syst(E_0^{\oplus n},\gamma) 
\times_{M_H(\gamma)} M_H(\gamma)_i \to M_H(\gamma)_i$
is $Gr(i,n)$-bundle.

%Then the Zariski closure $\overline{M_H(\gamma)_i}$
%of $M_H(\gamma)_i$ is $\cup_{j \geq i}M_H(\gamma)_j$
%if $i \geq \chi(\gamma_0,\gamma)$.

\begin{lem}\cite[Lem. 2.1]{Y:5}\label{lem:basic}
For $E \in M_H(\gamma)$ and $V \subset \Hom(E_0,E)$,
\begin{enumerate}
\item 
$ev:V \otimes E_0 \to E$ is injective and $\coker(ev)$ is stable.
\item
$ev:V \otimes E_0 \to E$ is surjective in codimension 1
and $\ker(ev)$ is stable.
\end{enumerate}
\end{lem} 

\begin{lem}\label{lem:D}
If $ev:V \otimes E_0 \to E$ is surjective in codimension 1,
then 
\begin{enumerate}
\item
$D(E):={\cal E}xt^1(V \otimes E_0 \to E,{\cal O}_X)$
is a stable sheaf of $\deg_{E_0^{\vee}}D(E)=1$.
\item
$\Ext^1(E_0,E)=0$.
\end{enumerate}
In particular $\chi(\gamma_0,\gamma) \geq n$.
\end{lem}

\begin{proof}
We have an exact sequence
\begin{equation}
{\cal E}xt^1(\im(ev) \to E,{\cal O}_X) \to
{\cal E}xt^1(V \otimes E_0 \to E,{\cal O}_X) \to 
{\cal H}om(\ker(ev),{\cal O}_X) \to
{\cal E}xt^2(\im(ev) \to E,{\cal O}_X).
\end{equation}
By Lemma \ref{lem:basic}, $\ker (ev)$ is stable and 
$\coker(ev)$ is of $0$-dimensional.
Then ${\cal E}xt^1(\im(ev) \to E,{\cal O}_X) \cong
{\cal E}xt^1(\coker(ev),{\cal O}_X)=0$ and
${\cal E}xt^2(\im(ev) \to E,{\cal O}_X) \cong
{\cal E}xt^2(\coker(ev),{\cal O}_X)$ is of $0$-dimensional.
Hence $D(X)$ is stable.

We next show that $\Ext^1(E_0,E)=0$.
Since $\ker(ev)$ is stable,
we get 
\begin{equation}
\Ext^2(E_0,\ker(ev))=\Hom(\ker(ev),E_0(K_X))^{\vee}=0.
\end{equation}
Combining the fact $\Ext^1(E_0,E_0)=0$, we see that
$\Ext^1(E_0,\im(ev))=0$.
Since $\Ext^1(E_0,\coker(ev))=0$, we get $\Ext^1(E_0,E)=0$.
\end{proof}

\begin{prop}
$\Syst(E_0^{\oplus n},\gamma)$ is smooth and
$\dim \Syst(E_0^{\oplus n},\gamma)=
\dim M_H(\gamma)-n(n-\chi(\gamma_0,\gamma))$.
\end{prop}

\begin{proof}
Let $(E,V) \in \Syst(E_0^{\oplus n},\gamma)$ be a 
$E_0$-twisted coherent system.
Since $V \subset \Hom(E_0,E)$, we have a homomorphism
\begin{equation}
\Hom(V \otimes E_0,V \otimes E_0) \to
\Hom(V \otimes E_0,E) \to \Ext^1(V \otimes E_0 \to E,E).
\end{equation}
Then the cokernel is the Zariski tangent space of
$\Syst(E_0^{\oplus n},\gamma)$ and
the obstruction space is $\Ext^2(V \otimes E_0 \to E,E)$.
If $\rk(\gamma-n\gamma_0) \geq 0$, then $\Ext^2(V \otimes E_0 
\overset{ev}{\to} E,E) \cong \Ext^2(\coker(ev),E)=0$.
If $\rk(\gamma-n\gamma_0) < 0$, then by using Lemma \ref{lem:D} and
an exact sequence
\begin{equation}
 \Ext^1(V \otimes E_0,E) \to \Ext^2(V \otimes E_0 \to E,E) \to \Ext^2(E,E),
\end{equation}
we see that $\Ext^2(V \otimes E_0 \to E,E)=0$.
Hence $\Syst(E_0^{\oplus n},\gamma)$ is smooth.
Then we see that
\begin{equation}
\begin{split}
\dim \Syst(E_0^{\oplus n},\gamma)&=
\dim \Ext^1(V \otimes E_0 \to E,E)-\dim PGL(V)\\
&=-\chi(E,E)+n\chi(E_0,E)-n^2\\
&=\dim M_H(\gamma)-n(n-\chi(\gamma_0,\gamma)).
\end{split}
\end{equation}
\end{proof}

\begin{prop}\label{prop:Gr}
We set $m:=-\chi(\gamma,\gamma_0)$.
\begin{enumerate}
\item
If $\rk \gamma \geq n \rk \gamma_0$, then
$\Syst(E_0^{\oplus n},\gamma)$ is a $Gr(m+n,n)$-bundle over 
$M_H(\gamma-n\gamma_0)$. 
\item
If $\rk \gamma < n \rk \gamma_0$, then
$\Syst(E_0^{\oplus n},\gamma) \cong 
\Syst((E_0^{\vee})^{\oplus n},n\gamma_0^{\vee}-\gamma^{\vee})$. 
In particular
$\Syst(E_0^{\oplus n},\gamma)$ 
is a $Gr(m+n,n)$-bundle over
$M_H(n\gamma_0^{\vee}-\gamma^{\vee})$. 
\end{enumerate}
\end{prop}

\begin{proof}
We first assume that $\rk \gamma \geq n \rk \gamma_0$.
For $(E,V) \in \Syst(E_0^{\oplus n},\gamma)$,
Lemma \ref{lem:basic} implies that
$ev:V \otimes E_0 \to E$ is injective and $\coker(ev)$ is stable.
Thus we have a morphism $\pi_n:\Syst(E_0^{\oplus n},\gamma) \to
M_H(\gamma-n \gamma_0)$.
Conversely for $G \in M_H(\gamma-n\gamma_0)$ and
an $n$-dimensional subspace $U$ of $\Ext^1(G,E_0)$,
we have an extension
\begin{equation}
0 \to U^{\vee} \otimes E_0 \to E \to G \to 0
\end{equation}  
whose extension corresponds to the inclusion 
$U \hookrightarrow \Ext^1(G,E_0)$.
Then $E$ is stable. 
Since
$\dim \Ext^1(G,E_0)=-\chi(G,E_0)=-\chi(\gamma-n\gamma_0,\gamma_0)$
and there is a universal family,
we see that $\pi_n$ is a (Zariski locally trivial) $Gr(m+n,n)$-bundle.
Therefore we get our claim.

We next treat the second case.
For $(E,V) \in \Syst(E_0^{\oplus n},\gamma)$,
$D(E):={\cal E}xt^1(V \otimes E_0 \to E,{\cal O}_X)$
fits in an exact sequence
\begin{equation}
0 \to E^{\vee} \to (V \otimes E_0)^{\vee} \to
D(E) \to {\cal E}xt^1(E,{\cal O}_X) \to 0.
\end{equation}
Hence $(V \otimes E_0)^{\vee} \to D(E)$ defines a point of
$\Syst((E_0^{\vee})^{\oplus n},n\gamma_0^{\vee}-\gamma^{\vee})$.
Thus we get a morphism
$\psi:\Syst(E_0^{\oplus n},\gamma) \to
\Syst((E_0^{\vee})^{\oplus n},n\gamma_0^{\vee}-\gamma^{\vee})$.
Conversely for $(F,U) \in 
\Syst((E_0^{\vee})^{\oplus n},n\gamma_0^{\vee}-\gamma^{\vee})$,
we get a homomorphism
$U^{\vee} \otimes E_0 \to 
{\cal E}xt^1(U \otimes E_0^{\vee} \to F,{\cal O}_X)$. 
It gives the inverse of $\psi$ (for more details, see
\cite[Prop. 5.128]{K-Y:1}).
\end{proof}

\begin{lem}
\begin{enumerate}
\item
If $\rk(\gamma-\chi(\gamma_0,\gamma)\gamma_0) \geq 0$, then
$\overline{M_H(\gamma)_i} =\emptyset$
for $\rk(\gamma-i\gamma_0)<0$.
\item
If $\rk(\gamma-\chi(\gamma_0,\gamma)\gamma_0)<0$, then
$M_H(\gamma)_{\chi(\gamma_0,\gamma)}=M_H(\gamma)$.
\end{enumerate}
\end{lem}

\begin{proof}
If $\dim(E_0,E)=i$ with $\rk(\gamma-i\gamma_0)<0$,
then Lemma \ref{lem:D} implies that $\chi(\gamma_0,\gamma) \geq i$.
Hence $\rk(\gamma-\chi(\gamma_0,\gamma)\gamma_0)<0$.
By Lemma \ref{lem:D},
 $\Ext^1(E_0,E)=0$ for all $E \in M_H(\gamma)$.
Hence $M_H(\gamma)_{\chi(\gamma_0,\gamma)}=M_H(\gamma)$.
\end{proof} 
By using Proposition \ref{prop:Gr}, we get the following theorem.

\begin{thm}
We set $\zeta:=\gamma(L_{e_0}(e))=\gamma-\chi(\gamma,\gamma_0) \gamma_0$
and $s:=-(K_X,c_1(e_0^{\vee} \otimes e))$.
Assume that $n:=-\chi(\gamma,\gamma_0)>0$.
Then
$M_H(\gamma) \cong \Syst(E_0^{\oplus n},\zeta)$
and we get a morphism $\lambda_{\gamma_0,\gamma}:M_H(\gamma)
\to M_H(\zeta)$
by sending $E$ to a universal extension
\begin{equation}
0 \to E_0 \otimes \Ext^1(E,E_0)^{\vee} \to 
\lambda_{\gamma_0,\gamma}(E) \to E \to 0.
\end{equation}
Hence we have a stratification
\begin{equation}
M_H(\gamma)= \coprod_{i \geq s}
\lambda_{\gamma_0,\gamma}^{-1}(M_H(\zeta)_i)
\end{equation}
such that $\lambda_{\gamma_0,\gamma}^{-1}(M_H(\zeta)_i) \to M_H(\zeta)_i$
is a $Gr(i,n)$-bundle.
In particular, $M_H(\gamma)_0 \to M_H(\zeta)_{n}$ is an isomorphism
for $n \geq s$.
\end{thm}

\begin{cor}
If $0>\chi(e_0,e)=-k\geq -s$, then
$M_H(\gamma(e))\to M_H(\gamma(L_{e_0}(e)))$ 
is birationally $Gr(s,k)$-bundle.
In particular, if $\chi(e_0,e)=-s$, then 
$M_H(\gamma(e)) \to M_H(\gamma(L_{e_0}(e)))$ is a birational map.
\end{cor}

\vspace{1pc}

\begin{ex}
Assume that $(X,H)=({\Bbb P}^1 \times {\Bbb P}^1,
{\cal O}_{{\Bbb P}^1 \times {\Bbb P}^1}(1,n))$, $n>0$.
We set $L:={\cal O}_{{\Bbb P}^1 \times {\Bbb P}^1}(-1,n+1)$.
Then $(L,H)=1$, $s=(L,-K_X)=2n$ and $\chi(L)=0$.
Hence $M_H(1+r,L,r) \cong Gr(2n,r)$.
\end{ex}

\subsection{Virtual Hodge polynomial}
We set $a:=-\chi(\gamma,\gamma_0)$.
Assume that $\rk(\gamma-\chi(\gamma_0,\gamma)\gamma_0) \geq 0$.
We shall consider vitrual Hodge polynomial of
$M_H(\gamma+k \gamma_0)_i$.
For an algebraic set $Z$, 
\begin{equation}
e(Z):=\sum_{p,q}(-1)^{p+q}h^{p,q}(Z)x^p y^q
\end{equation}
is the virtual Hodge polynomial of $Z$
(cf. \cite{D-K:1}).
We set $t:=xy$.
Then
\begin{equation}
\begin{split}
e(M_H(\gamma+k \gamma_0)_j)&=
e(Gr(a+j-k,j))e(M_H(\gamma+(k-j) \gamma_0)_0)\\
&=\frac{[a+j-k]!}{[a-k]![j]!}e(M_H(\gamma+(k-j) \gamma_0)_0),
\end{split}
\end{equation}
where 
\begin{equation}
[n]:=\frac{t^n-1}{t-1}, \quad [n]!:=[n][n-1]\cdots[1].
\end{equation}
By summing up all $e(M_H(\gamma+k \gamma_0)_k)$,
we get
\begin{equation}
\sum_k [a-k]! e(M_H(\gamma+k \gamma_0))y^k
=\left(\sum_j \frac{1}{[j]!}y^j \right) 
\left(\sum_l [a-l]! e(M_H(\gamma+l \gamma_0)_0)y^l \right).
\end{equation}
Since 
\begin{equation}
\left(\sum_j \frac{1}{[j]!}y^j \right)^{-1}=
\sum_j \frac{(-1)^jt^{j(j-1)/2}}{[j]!}y^j, 
\end{equation} 
we get that
\begin{lem}
If $\rk(\gamma-\chi(\gamma_0,\gamma)\gamma_0) \geq 0$, then
\begin{equation}
e(M_H(\gamma+l \gamma_0)_0)=\sum_{j \geq 0} (-1)^{j} t^{j(j-1)/2} 
\frac{[a+j-l]!}{[a-l]![j]!}e(M_H(\gamma+(l-j) \gamma_0)).
\end{equation}
In particular
\begin{equation}
e(M_H(\gamma+k \gamma_0)_i)=\sum_{j \geq 0} (-1)^{j} t^{j(j-1)/2} 
\frac{[a-k+i+j]!}{[a-k]![i]![j]!}e(M_H(\gamma+(k-i-j) \gamma_0)).
\end{equation}
Since $M_H(\gamma+l \gamma_0)_0 = \emptyset$ for $a-s<l \leq a$,
we also get the following relations:
\begin{equation}
\sum_{j \geq 0} (-1)^{j} t^{j(j-1)/2} 
\frac{[a+j-l]!}{[a-l]![j]!}e(M_H(\gamma+(l-j) \gamma_0))=0
\end{equation}
for $a-s<l \leq a$.
\end{lem}

\subsection{Examples on ${\Bbb P}^2$}

From now on, we assume that $X$ is ${\Bbb P}^2$. 
Then $s=-(K_X,{\cal O}_X(1))=3$.
Hence we get the following relations:
\begin{equation}\label{eq:relation2}
\begin{split}
&\sum_{j \geq 0} (-1)^{j} t^{j(j-1)/2} 
e(M_H(\gamma+(a-j) \gamma_0))=0,\\
&\sum_{j \geq 0} (-1)^{j} t^{j(j-1)/2} 
[j+1]e(M_H(\gamma+(a-1-j) \gamma_0))=0,\\
&\sum_{j \geq 0} (-1)^{j} t^{j(j-1)/2} 
\frac{[j+2][j+1]}{[2]!}e(M_H(\gamma+(a-2-j) \gamma_0))=0.
\end{split}
\end{equation}

By a simple calculation, we get
\begin{prop}\label{prop:relation}
\begin{equation}\label{eq:relation}
\begin{split}
e(M_H(\gamma+(a-2) \gamma_0)) & =\sum_{j \geq 0} (-1)^{j} t^{(j+1)j/2} 
\frac{[j+3][j+2]}{[2]!}e(M_H(\gamma+(a-3-j) \gamma_0)),\\
e(M_H(\gamma+(a-1) \gamma_0))& =\sum_{j \geq 0} (-1)^{j} 
t^{(j+1)j/2}[j+3][j+1]e(M_H(\gamma+(a-3-j) \gamma_0)),\\
e(M_H(\gamma+a \gamma_0))&=\sum_{j \geq 0} (-1)^{j} t^{(j+1)j/2} 
\frac{[j+2][j+1]}{[2]!}e(M_H(\gamma+(a-3-j) \gamma_0)).
\end{split}
\end{equation}
\end{prop}

Assume that $E_0:={\cal O}_X$.
We set $\gamma:=\gamma({\cal O}_X(1))$.
Then
\begin{equation}
\begin{split}
M_H(\gamma-a \omega-\gamma_0)&=\{{\cal O}_l(1-a) |
\text{ $l$ is a line on ${\Bbb P}^2$} \}\\
& \cong {\Bbb P}^2.
\end{split}
\end{equation}
Hence $M_H(\gamma-a \omega)_1$, $a \geq 2$ is a ${\Bbb P}^a$-bundle over
${\Bbb P}^2$. 
By the morphism $M_H(\gamma-a \omega) \to M_H(\gamma-a\omega+a\gamma_0)$,
the fibers of $M_H(\gamma-a \omega)_1 \to {\Bbb P}^2$ are contracted.

\vspace{1pc}

\begin{ex}
If $a=2$, then $M_H(\gamma-2\omega+2\gamma_0) \cong M_H(\gamma^2-\gamma_0) \cong
{\Bbb P^2}$.
That is, $E \in M_H(\gamma-2\omega+2\gamma_0)$ fits in a universal extension
\begin{equation}
0 \to {\cal O}_X^{\oplus 3} \to E \to {\cal O}_l(-1) \to 0.
\end{equation}
Moreover we see that
$M_H(\gamma-2\omega+i\gamma_0)$, $i=0,1$ are ${\Bbb P}^2$-bundle over
$M_H(\gamma-2\omega+2\gamma_0) \cong {\Bbb P}^2$.
\end{ex}

\vspace{1pc}

\begin{ex}
If $a=3$, then $M_H(\gamma-3\omega) \to M_H(\gamma-3\omega+3\gamma_0)$ 
is the blow-up of
$M_H(\gamma-3\omega+3\gamma_0)_4 \cong M_H(\gamma-3\omega-\gamma_0)$.
This was obtained by Drezet \cite[IV]{D:3}.
\end{ex}

By \cite{E-S:1} and \cite{Y:1},
we know $e(M_H(r,H,\chi))$ for $r=1,2$.
By using Proposition \ref{prop:relation}, we get the following:

\begin{equation}
\begin{split}
e(M_H(1,H,0))&=1+2t+5t^2+6t^3+5t^4+2t^5+t^6,\\
e(M_H(2,H,1))&=1+2t+6t^2+9t^3+12t^4+9t^5+6t^6+2t^7+t^8,\\
e(M_H(3,H,2))&=1+2t+5t^2+8t^3+10t^4+8t^5+5t^6+2t^7+t^8,\\
e(M_H(4,H,3))&=1+t+3t^2+3t^3+3t^4+t^5+t^6.
\end{split}
\end{equation}

\begin{equation}
\begin{split}
e(M_H(1,H,-1))&=1+2t+6t^2+10t^3+13t^4+10t^5+6t^6+2t^7+t^8,\\
e(M_H(2,H,0))&=1+2t+6t^2+13t^3+24t^4+35t^5+41t^6+
 35t^7+24t^8+13t^9+6t^{10}+2t^{11}+t^{12},\\
e(M_H(3,H,1))&=1+2t+6t^2+12t^3+24t^4+38t^5+54t^6+
 59t^7+54t^8+38t^9+24t^{10}+12t^{11}+6t^{12}+2t^{13}+t^{14},\\
e(M_H(4,H,2))&=1+2t+5t^2+10t^3+18t^4+28t^5+38t^6+
 42t^7+38t^8+28t^9+18t^{10}+10t^{11}+5t^{12}+2t^{13}+t^{14},\\
e(M_H(5,H,3))&=1+t+3t^2+5t^3+8t^4+10t^5+12t^6+
 10t^7+8t^8+5t^9+3t^{10}+t^{11}+t^{12}.
\end{split}
\end{equation}

If $E_0:=\Omega_X(1)$, then $\deg_{E_0}({\cal O}_X)=1$.
We set $\gamma=\gamma({\cal O}_X)$. Then
\begin{itemize}
\item
$M_H(\gamma-a \omega) \to M_H(\gamma-a \omega+2a \gamma_0)$ is a closed 
immersion for $a \geq 2$.
\item
If $a=2$, then $M_H(\gamma-2\omega+\gamma_0) \to 
M_H(\gamma-2\omega+4 \gamma_0)$
is the blow-up along $M_H(\gamma-2\omega)$.
\end{itemize}
Here we remark that
Drezet showed that $M_H(\gamma-2\omega+4\gamma_0)=M_H(9,-4H,-1)
 \cong Gr(6,2)$
(see \cite[Appendice]{D:1}).
Since $e(M_H(1,0,-1))=1+2t+3t^3+2t^3+t^4$ and
$e(M_H(3,-H,-1))=e(M_H(3,H,2))$,
Proposition \ref{prop:relation} implies that
\begin{equation}
\begin{split}
e(M_H(3,-H,-1))&=1+2t+5t^2+8t^3+10t^4+8t^5+5t^6+2t^7+t^8,\\
e(M_H(5,-2H,-1))&=1+2t+5t^2+8t^3+13t^4+14t^5+13t^6+8t^7+5t^8+2t^9+t^{10},\\
e(M_H(7,-3H,-1))&=1+2t+4t^2+6t^3+9t^4+10t^5+9t^6+6t^7+4t^8+2t^9+t^{10},\\
e(M_H(9,-4H,-1))&=1+t+2t^2+2t^3+3t^4+2t^5+2t^6+t^7+t^8
(=e(Gr(6,2))).
\end{split}
\end{equation}

\subsubsection{Line bundles on $M_H(\gamma)$}

Let $p_{M_H(\gamma(e))}:M_H(\gamma(e)) \times X \to
M_H(\gamma(e))$ and $q:M_H(\gamma(e)) \times X \to X$ be projections,
and let
${\cal E}$ be a universal family on $M_H(\gamma(e)) \times X$.
We define a homomorphism
 $\theta_e:e^{\perp} \to \Pic(M_H(\gamma(e)))$ by
\begin{equation}
\theta_e(x):=\det p_{M_H(\gamma(e))!}({\cal E}^{\vee} \otimes q^*(x)),
\end{equation}
where $e^{\perp}:=\{x \in K(X)|\chi(e,x)=0 \}$.
The following is a special case of Drezet's results.
\begin{thm}\cite{D:2}
Assume that $\dim M_H(\gamma(e))=1-\chi(e,e)>0$.
Then $\theta_e$ is surjective and
\begin{enumerate}
\item 
$\theta_e$ is an isomorphism, if $\chi(e,e)<0$,
\item
$\ker \theta_e={\Bbb Z}e_0$, if $\chi(e,e_0)=0$.
\end{enumerate}
\end{thm}
We set $\tilde{e}:=L_{e_0}(e)$.
By a simple calculation, we see that the following diagram
is commutative:
\begin{equation}
\begin{CD}
e^{\perp} @<{R_{e_0}}<< \tilde{e}^{\perp}/e_0\\
@V{\theta_e}VV @VV{\theta_{\tilde{e}}}V \\
\Pic(M_H(\gamma(e))) @<<{\lambda_{\gamma(e_0),\gamma(e)}^*}< 
\Pic(M_H(\gamma(\tilde{e})))
\end{CD}
\end{equation}
We set $\alpha_e:=-(\rk e) {\cal O}_H+\chi(e,{\cal O}_H){\Bbb C}_P$.
Then it gives a map to the Uhlenbeck compactification
\cite{Li:1}.
$\beta_e:=R_{e_0}(\alpha_{\tilde{e}})$
gives the map $\lambda_{\gamma(e_0),\gamma(e)}:
M_H(\gamma(e)) \to M_H(\gamma(\tilde{e}))$.
\begin{itemize}
\item
If $E_0={\cal O}_X$, $\rk e>0$ and $\chi(e,e_0)<0$, 
then the nef. cone of $M_H(\gamma(e))$ is generated by 
$\alpha_e$ and $\beta_e$. 
\end{itemize}
This is a generalization of \cite{S:1}. 

For $\gamma:=(3,H,5-a)$, 
we set $\gamma_0:=(1,0,1)$,
$\gamma_1:=\gamma(\Omega_X(1))=(2,-H,0)$,
$\delta:=\gamma+a\gamma_0$ and
$\eta:=\gamma^{\vee}+(2a-3)\gamma_1$.
Then we get the following diagram:

\begin{equation}
\begin{matrix}
           && {M}_H(\gamma)&&\leftarrow \cdots \rightarrow &&
    {M}_{H}(\gamma^{\vee})&&&\cr
          &\llap{$\lambda_{\gamma_0,\gamma}$}
  \swarrow&&\searrow\rlap{}&&  
     \llap{}\swarrow 
    &&\searrow\rlap{$\lambda_{\gamma_1,\gamma^{\vee}}$}\cr
          M_H(\delta) && 
   && N_H(\gamma)&&&&
{M}_H(\eta)\cr    
\end{matrix}
\end{equation} 
where $N_H(\gamma)$ is the Uhlenbeck compactification of
$M_H(\gamma)^{\mu\text{-}s,loc}$.
$M_H(\gamma^{\vee})$ contains ${\Bbb P}^{2a-3}$-bundle over
$M_H(1,0,2-a)$ and 
$\lambda_{\gamma_0,\gamma^{\vee}}$ contracts the fibers.
${\lambda_{\gamma_0,\gamma}}_{|M_H(\gamma)_i}$ is a $Gr(a-2+i,a-2)$-bundle
over $M_H(\delta)_{a-2+i} \cong M_H(\gamma-i\gamma_0)_0$.
Then it is easy to see that
$M_H(3,H,5-a) \not \cong M_H(3,-H,2-a)$.

\vspace{1pc}

{\it Acknowledgement.}
A starting point of this note is \cite{E-S:1}.
I would like to thank M. Maruyama for giving me
a preprint version of \cite{E-S:1}.

\end{document}